\documentclass[12pt]{article}
\usepackage{bbm}
\usepackage{graphicx}
\graphicspath{ {./images/} }
\usepackage{wrapfig}
\usepackage[mathscr]{euscript}
\usepackage{subcaption}
\usepackage{hyperref}
\usepackage{epstopdf}
\usepackage{enumerate}
\usepackage{amsmath,amsfonts,amssymb,amsthm,epsfig,epstopdf,titling,url,array}
\usepackage{color}
\usepackage{framed}
\usepackage[utf8]{inputenc}
\usepackage[english]{babel}
\usepackage{xparse}
\usepackage{authblk}
\usepackage[margin=2.5cm]{geometry}
\usepackage{authblk}
\newtheorem{thm}{Theorem}[section]

\newtheorem{lem}{Lemma}[section]

 \newtheorem{rem}{Remark}[section]

\date{}
 
\begin{document}
\title{On the entire solutions of Schr\"{o}der iteration functions of first and second kind}


  \author[1]{Subhasis Ghora
  	\footnote{subhasisghora@iiitdm.ac.in}}
  \author[2]{Gang Liu
  	\footnote{liugangmath@sina.cn}}
  \author[3]{Soumen Pal \footnote{soumen.pal.new@gmail.com}}
  \author[3]{Saminathan Ponnusamy  \footnote{samy@iitm.ac.in}}
  \affil[1]{Department of Sciences and Humanities, Indian Institute of Information Technology Design and Manufacturing Kancheepuram, Chennai, India}
  \affil[2]{College of Mathematics and Statistics, Hengyang Normal University,	Hengyang, China}
  \affil[3]{Department of Mathematics,
  	Indian Institute of Technology Madras, Chennai, India}
	\date{}
\maketitle
\begin{abstract}
Let $S_{f,n}$ and $K_{f,n}$ denote Schr\"oder's methods of the first and second kind with given order $ n $ ($\geq 2 $)  applied to an entire function $ f $. This article investigates the existence of the entire solutions $f$ whenever $ S_{f,n} $ and $ K_{f,n} $ are rational. Utilizing the Nevanlinna theory, we obtain the general form of entire solutions for such equations. It is proved that either of $ S_{f,n} $ and $ K_{f,n} $ is a non-identity  rational map if and only if   $f = p e^{q}$,  for some  polynomials $p$ and $q$.  Additionally, we provide one recurrence relation of $K_{f,n}$ depending only upon  $K_{f,n-1}$.
\end{abstract}
\textit{Keywords:}  Schr\"{o}der's methods; Exponential function; Nevanlinna characteristic function.\\
AMS Subject Classification: 30D35, 65H05
\section{Introduction}
Among the iterative methods for solving nonlinear equations, Schr\"oder's methods of first and second kind form two classical and remarkably broad families. Introduced by Schr\"oder in 1870 \cite{sch1} (see also the translation \cite{sch2}), these methods provide systematic higher-order generalizations of the Newton-type schemes. They have since appeared in several forms in numerical analysis, approximation theory, and the study of functional equations. Beyond their role in root-finding, Schr\"oder-type constructions are closely related to the linearization of dynamical systems, acceleration of convergence, inverse function approximation, and iterative functional equations. Consequently, they arise naturally in several areas including complex dynamics, symbolic computation, optimization, and computational mathematics.

For an entire function, Schr\"oder's method of first kind of order $n$ ($\geq 2$), denoted by $S_{f,n}$, is defined by
\begin{equation}
S_{f,n}(z)=z+\sum_{k=1}^{n-1}c_k(z)(-f(z))^k,
\label{sch-def1}
\end{equation}
where
\[
c_k(z)=\frac{1}{k!}\left(\frac{1}{f'(z)}\frac{d}{dz}\right)^{k-1}\frac{1}{f'(z)}.
\]
It can also be expressed recursively as follows. Define
\begin{equation}
s_{f,1}(z)=1,\qquad
s_{f,m+1}(z)=f'(z)s_{f,m}'(z)-(2m-1)f''(z)s_{f,m}(z),
\label{rr}
\end{equation}
and then
\begin{equation}
S_{f,n}(z)=z+\sum_{k=1}^{n-1}\frac{(-1)^k}{k!}\frac{s_{f,k}(z)}{(f'(z))^{2k-1}}(f(z))^k.
\label{sch-def2}
\end{equation}
Schr\"oder's method of second kind of order $n$, denoted by $K_{f,n}$, is given by
\begin{equation}
K_{f,n}(z)=z+(n-1)\frac{\left(\frac{1}{f(z)}\right)^{(n-2)}}{\left(\frac{1}{f(z)}\right)^{(n-1)}}
\label{sch2-def11}
\end{equation}
which can also be represented using a recurrence relation. Define
\begin{equation}
\nonumber
l_{f,1}(z)=1, ~\mbox{and~}l_{f,k+1}(z)=f(z)l_{f,k}'(z)-kf'(z)l_{f,k}(z).
\label{rr2}\end{equation}
An induction using this recursion gives
$$\left(\frac1f\right)^{(k-1)}=\frac{l_{f,k}}{f^k}\qquad(k\ge1).$$
Then the formula for $K_{f,n}$ can be written as
\begin{equation}
\nonumber
K_{f,n}(z)=z+(n-1)\frac{l_{f,n-1}(z)f(z)}{l_{f,n-1}(z)}.
\label{sch2-def2}\end{equation}
These two families contain several well-known methods as special cases. In particular,
\[
S_{f,2}(z)=K_{f,2}(z)=N_f(z)=z-\frac{f(z)}{f'(z)},
\]
the classical Newton map of $f$. For $n=3$, one obtains Chebyshev's method
\[
C_f(z)=z-\left(1+\frac{1}{2}\frac{f(z)f''(z)}{(f'(z))^2}\right)\frac{f(z)}{f'(z)}
\]
from the first kind, and Halley's method
\[
H_f(z)=z-\frac{2f'(z)f(z)}{2(f'(z))^2-f(z)f''(z)}
\]
from the second kind.  The dynamics of these methods have been studied extensively in the literature; see, for example, \cite{hen,DS2022,LPP2026,Nayak-Pal2022}, and the references therein.

These two families have been extensively studied from the viewpoints of numerical efficiency and local convergence; see for instance \cite{coll,sir,pet,xin}. Their global behavior as dynamical systems, however, remains much less understood, especially for transcendental entire functions and for higher-order members of the two Schr\"oder families.

From the viewpoint of complex dynamics, an important issue is whether the map associated with a root-finding method is \emph{rational} or \emph{transcendental meromorphic}. This distinction plays a fundamental role in determining the nature of the iteration. If the associated map is rational, then its behaviour is governed by the classical theory of rational iteration. In particular, Sullivan's no wandering domain theorem \cite[Theorem 8.1.2]{Beardon_book} ensures that wandering domains do not occur in this setting. On the other hand, if the associated map is transcendental meromorphic, then the dynamics can be significantly more complicated. In particular, Baker domains \cite{Beardon_book} and wandering domains \cite{Beardon_book} may arise, and the presence of such invariant or non-recurrent regions can obstruct the root-finding process. As a property of root-finding methods, every root of $f$ is an attracting fixed point of both $S_{f,n}$ and $K_{f,n}$, and so one hopes that a typical orbit converges to a root. But if an orbit starts in a Baker domain or a wandering domain, then it may never approach any attracting fixed point at all, thereby defeating the very purpose of the method. Thus, determining when these higher-order Schr\"oder maps are rational, is not merely a structural question: it is directly relevant to the global reliability of the corresponding root-finding methods.

To place this in this context for Newton's method, we recall that the associated map $N_f$ is rational if and only if
$$f(z)=p(z)e^{q(z)}$$
for two polynomials $p$ and $q$; see \cite[Proposition 2.11]{RC2007}. More recently, analogous characterizations have been obtained for Halley's method \cite[Theorem 3.6]{CMHD2022} and for Chebyshev's method \cite[Theorem A]{pooja}. These results reveal a striking and rigid phenomenon: among all entire functions, only those of the form $p e^q$ give rise to rational root-finding maps. At the same time, the proofs for Newton, Halley, and Chebyshev rely on rather different ideas. Newton's case is comparatively elementary, Halley's method leads to a delicate first-order differential equation, and the proof for Chebyshev's method requires tools from Nevanlinna theory. This naturally raises the following broad question.

\medskip
\noindent
\textit{For the higher-order Schr\"oder families of first and second kind, which entire functions give rise to rational maps?}
\medskip

This question is natural for several reasons. Firstly, Schr\"{o}der's methods provide a unified framework that includes Newton's method, Chebyshev's method, Halley's method, and many other higher-order iterative methods. Therefore, a characterization at this level gives a common structural understanding of a broad class of root-finding maps. Secondly, unlike the low-order cases, the higher-order formulas involve increasingly complicated differential expressions, and so it is not immediately clear whether the same structural behaviour persists. Thirdly, from a dynamical viewpoint, such a characterization tells us precisely when the transcendental complications discussed above are absent.

Our first main result shows that the same class of functions, namely, those of the form $p e^q$, characterizes all entire functions whose Schr\"oder maps of first kind are rational.

\begin{thm}\label{first}
	Let $f:\mathbb{C}\to\mathbb{C}$ be a non-constant entire function. Then Schr\"oder's method of first kind of order $n$ applied to $f$ is a non-identity rational map if and only if
$$f(z)=p(z)e^{q(z)}$$
	for some polynomials $p\,(\not\equiv 0)$ and $q$.
\end{thm}
For the second kind, the situation is more delicate. A key step in our analysis is the direct logarithmic-derivative argument.   A new recurrence relation connecting two successive orders:
$$
K_{f,n+1}(z)=z-n\frac{z-K_{f,n}(z)}{\,n-K'_{f,n}(z)\,},\qquad n\geq 2.
$$
This relation appears to be of independent interest, since it expresses the $(n+1)$-th member of the family entirely in terms of the $n$-th one and its derivative.  We state  our second main result.
\begin{thm}\label{sec}
	Let $f:\mathbb{C}\to\mathbb{C}$ be a non-constant entire function. Then Schr\"oder's method of second kind of order $n$ applied to $f$ is rational if and only if
$$	f(z)=p(z)e^{q(z)}$$
	for some polynomials $p\,(\not\equiv 0)$ and $q$.
\end{thm}
Theorems~\ref{first} and \ref{sec} extend, in a unified way, the well-known rationality characterizations for Newton's, Halley's, and Chebyshev's methods \cite{CMHD2022,Nayak-Pal2022,RC2007} to the full higher-order Schr\"oder families. In particular, they show that the class of entire functions of the form $p e^q$ is not a feature of only a few isolated methods, but is in fact the natural class governing rationality throughout these two fundamental hierarchies of iterative methods.

The rest of the paper is organized as follows. In Section~\ref{2_Prel}, we discuss preliminaries and notation, including the characteristic function and some related results. In Section~\ref{3_Proof}, we prove the main results, namely, Theorems~\ref{first} and \ref{sec}.
\section{Preliminaries and notation}\label{2_Prel}
For two functions $f$ and $g$ from the positive real line into itself, we say $f(x)=o(g(x))$ as $x \to \infty$ if $\lim\limits_{x \to \infty} \frac{f(x)}{g(x)}=0$. We say $f(x)=\mathcal{O}(g(x))$ as $x \to \infty$ if  $ \frac{f(x)}{g(x)} \leq K$ for some $K>0$ and for all sufficiently large positive $x$. Note that $f(x)=o(g(x))$ implies $f(x)=\mathcal{O}(g(x))$, but not conversely.

 We recall some standard facts from Nevanlinna theory that will be used later in the paper. Let $f$ be a non-constant meromorphic function in the disc $|z|\leq R$ ($0<R<\infty$). For $0<r<R$, we define
 $$m(r,f)=\frac{1}{2\pi}\int_0^{2\pi}\log^+|f(re^{i\theta})|d\theta,$$  where $\log^+(x)=\max\{\log x,0\}$ and
 $$N(r,f)=\int_0^r\frac{n(t,f)-n(0,f)}{t}dt+n(0,f)\log r,$$
where $n(t,f)$ denotes the number of poles of $f$ in the disc $|z|\leq t$ counting multiplicity. The characteristic function of $f$, denoted by $T(r,f)$, is defined as $$T(r,f)=m(r,f)+N(r,f).$$ It is well-known that if $f_1, f_2,\dots,f_n$ are meromorphic then for $r\geq 1$, we have the following inequalities: $$T(r,f_1f_2\dots f_n)\leq \sum_{i=1}^{n}T(r,f_i)$$  and $$T(r,\sum_{i=1}^{n}f_i)\leq \sum_{i=1}^{n}T(r,f_i)+\ln n.$$ 
Moreover, it follows from \cite[Proposition~2.1.11]{lai1993} that if $\alpha, \beta, \gamma, \delta \in \mathbb{C}$ and  $f$ is meromorphic such that $f\neq -\frac{\delta}{\gamma}$ then $$T(r, \frac{\alpha f+\beta}{\gamma f+\delta})=T(r,f)+\mathcal{O}(1).$$
\begin{rem}

\begin{enumerate}
	\item If $w$ is meromorphic and $j\ge1$, then
	\begin{equation}\label{log_derivative}
	m\left(r,\frac{w^{(j)}}{w}\right)
	=
	\mathcal{O}\left(\log r+\log^+T(r,w)\right)
	\end{equation}
	for all $r$ outside a set $E\subset(0,\infty)$ of finite linear
	measure.
	
	\item We shall also use the elementary fact that if $R$ is rational, then
	\begin{equation}\label{rational_proximity}
	m(r,R)=\mathcal{O}(\log r).
	\end{equation}
	\item If $w$ is entire, then
	\begin{equation}\label{max_modulus}
	\log M(r,w)
	\le
	3T(2r,w)+\mathcal{O}(1).
	\end{equation}
	
\end{enumerate}

\end{rem}

 Following two lemmas will be needed in the proof of our first main result.
\begin{lem}\label{lem2}
	Let $h$ be a non-constant entire function and $f(z)=e^{h(z)}$.
	Then we have $T(r,h')=o(T(r,f))$ and $T(r,h'')=o(T(r,f))$ as $r\rightarrow\infty$ and $r$ is outside a set of finite linear measure.
\end{lem}
\begin{proof}
Since $h$ is nonconstant, $f=e^{h}$ is a transcendental entire function. Hence, by \cite[Theorem~1.5]{yy2006},
\begin{equation}
	\lim_{r\to\infty}\frac{T(r,f)}{\log r}=\infty.
\label{c1}
\end{equation}
Clearly $T(r,h')=m(r, \frac{f'}{f})$. It follows from \cite[Lemma 1.4', p. 21]{yy2006},    that $T(r,h')$ is either $\mathcal{O}(\log r)$ or $\mathcal{O}(\log (rT(r,f)))$ outside of a set $E_0$ of finite linear measure. If $T(r,h')=\mathcal{O}(\log r)$, then using  Eqn.~(\ref{c1}), we have $T(r,h')=o(T(r,f))$ as $r\rightarrow\infty$. If $T(r,h')=\mathcal{O}(\log (rT(r,f)))$, then $T(r,h')=K\log (rT(r,f))$ for some $K>0$. Note that $$\lim\limits_{r\rightarrow\infty}\frac{\log (rT(r,f))}{T(r,f)}=\lim\limits_{r\rightarrow\infty}\frac{\log r+\log T(r,f)}{T(r,f)}=\lim\limits_{r\rightarrow\infty}\frac{\log r}{T(r,f)}+\lim\limits_{r\rightarrow\infty}\frac{\log T(r,f)}{T(r,f)}=0$$ outside the set $E_0$. This gives that $T(r,h')=o(T(r,f))$ as $r\rightarrow\infty$ outside of the set $E_0$. Using the general derivative estimate for meromorphic functions (cf. \cite[Theorem~1.47]{yy2006}), we have, $$T(r,h'') \le T(r,h') + \mathcal{O}(\log r + \log T(r,h')).$$ Since both $\log r$ and $\log T(r,h')$ are $o(T(r,f))$ by Eqn.~(\ref{c1}), it follows that $T(r,h'')=o(T(r,f))$ as $r\to\infty$ outside of the set $E_0$.
	\end{proof}
\begin{lem}\label{lem1}
\cite[Theorem~1.51]{yy2006}  For $n \geq 2$, suppose that $f_1, f_2, \dots, f_n $ are meromorphic functions and $g_1, g_2, \dots, g_n$ are entire functions satisfying the following conditions:
\begin{enumerate}
\item  $\sum\limits_{j=1}^n f_j(z) e^{g_j(z)} \equiv 0$, and $g_j(z)-g_k(z)$ are non-constants for all  $1 \leq j<k \leq n$.
\item For $1 \leq j \leq n $ and $ 1 \leq h<k \leq n$, we have $$ T\left(r, f_j\right)=o\left(T\left(r, e^{g_h-g_k}\right)\right) \quad(r \rightarrow \infty, r \notin E),$$
where $T(r,f)$ is the Nevanlinna characteristic function of $f$ and $E$ is a set with finite linear measure.
\end{enumerate}
Then $f_j(z) \equiv 0$ for all $j=1,2, \dots, n$.
\end{lem}
\section{Proofs of the main Theorems}\label{3_Proof}
In this section we present the proofs of the main results of this article.
\subsection{Schr\"{o}der iteration functions of first kind}
For an entire function $f$, consider Schr\"{o}der's method  of the first kind $S_{f,n} $ of order $n$ defined by Eqn.~(\ref{sch-def2}). A specific form of $s_{f,k}$ and its derivative is discussed in the following lemma for a particular type of entire function.
\begin{lem}
If $f(z)=p(z)e^{q(z)}$, where $p$ is a polynomial and $q$ is an entire function, then the functions $s_{f,k}$ and $s'_{f, k}$ can be expressed as $e^{(k-1)q(z)}P_{1,k}(z)$ and $e^{(k-1)q(z)}P_{2,k}(z)$, respectively, for some entire functions  $P_{1,k}$ and $P_{2,k}$.
Especially, the functions $P_{1,k}$ and $P_{2,k}$ are polynomials when $q$ is a polynomial.
	\label{sf1}
\end{lem}
\begin{proof}
We use the mathematical induction to prove the first part of the statement. It follows from Eqn.~(\ref{rr}) that the statement is trivially true for $k=1$. In this case, $P_{1,1}(z)=1$ and $P_{2,1}(z)=0$. For $k=2$, we have $$s_{f,2}(z)=e^{q(z)}P_{1,2}(z),$$ where $$P_{1,2}=-(p''+2p'q'+pq''+pq'^2)~ \text{and}~s'_{f,2}(z)=e^{q(z)}P_{2,2}(z)~ \text{with}~ P_{2,2}=P_{1,2}'+P_{1,2}q'.$$ Thus $P_{1,2}$ and $P_{2,2}$ are entire functions. Assuming that the statement is true for $k=m$, we have $$s_{f,m}(z)=e^{(m-1)q(z)}P_{1,m}(z)~~ \text{and}~~ s'_{f,m}(z)=e^{(m-1)q(z)}P_{2,m}(z)$$ for some entire functions $P_{1,m}$ and $P_{2,m}$. Now, for $k=m+1$, we have
\begin{align*}
s_{f,m+1}(z)&=f'(z)s_{f,m}'(z)-(2m-1)f''(z)s_{f,m}(z)\\
&=(p'+pq')(z)e^{q(z)}P_{2,m}(z)e^{(m-1)q(z)}
\\&~~~~-(2m-1)(p''+2p'q'+pq''+pq'^2)(z)e^{q(z)}e^{(m-1)q(z)}P_{1,m}(z)
\\&=e^{mq(z)}((p'+pq')P_{2,m}-(2m-1)(p''+2p'q'+pq''+pq'^2)P_{1,m})(z)
\\&=e^{mq(z)}P_{1,m+1}(z),
\end{align*}	
where $P_{1,m+1}=(p'+pq')P_{2,m}-(2m-1)(p''+2p'q'+pq''+pq'^2)P_{1,m}$.
Using above we have $$s'_{f,m+1}(z)=e^{mq(z)}P_{2,m+1}(z),$$ where $P_{2,m+1}=P'_{1,m+1}+mP_{1,m+1}q'$. Thus, the statement is true for $k=m+1$ whenever it is true for $k=m$.
\par
It follows from above proof that if $q$ is a polynomial then the functions $s_{f,k+1}$ and $s'_{f, k+1}$ can also be expressed as $e^{kq(z)}P_{1,k+1}(z)$ and $e^{kq(z)}P_{2,k+1}(z)$, respectively, where   $P_{1,k+1}$ and $P_{2,k+1}$ are polynomials.	
\end{proof}
The following lemma shows that, under our assumptions, the class of functions of the form $p e^q$ is preserved under differentiation.
\begin{lem}
Let $f$ be a non-constant entire function of the form $f(z)=p(z) e^{q(z)}$, where $p\,(\not\equiv 0)$ is a polynomial and $q$ is entire, and suppose that $S_{f,n}$ is a non-identity rational function. Then $f'$ is also of the form $p_1(z) e^{q_1(z)}$, where $p_1$ is a polynomial and $q_1$ is an entire  function.
\label{f'finitepole}
\end{lem}
%
%
%
%
%

\begin{proof}
	 It follows from Eqn.(\ref{rr}) and Eqn. (\ref{sch-def2}) that the functions $s_{f,k}$ are entire, and $S_{f,n}$ is the meromorphic function. The function $S_{f,n}$ can be written as 
	\begin{equation}\label{eq:S}
	S_{f,n}(z)-z=\sum_{k=1}^{n-1}\frac{(-1)^k}{k!}\,t_k(z),\quad
	\mbox{by considering~}t_k:=\frac{s_{f,k}\,f^{k}}{(f')^{2k-1}} .
	\end{equation}
	Considering $\psi:=f/f'$, so that $t_1=\psi$, we claim that
	\begin{equation}\label{eq:trec}
	t_{k+1}=\psi\,t_k'-k\,t_k\qquad(k\ge1).
	\end{equation}
	To proof it, note that on the open set where $ff'\neq0$, differentiating $t_k=s_{f,k}f^k(f')^{-(2k-1)}$ gives
	\[
	t_k'=\frac{s_{f,k}'f^{k}}{(f')^{2k-1}}+\frac{k\,s_{f,k}f^{k-1}}{(f')^{2k-2}}-\frac{(2k-1)\,s_{f,k}f^{k}f''}{(f')^{2k}} .
	\]
	Multiplying by $\psi=f/f'$ we have,
	\[
	\psi t_k'=\frac{f^{k+1}}{(f')^{2k+1}}\Big(f's_{f,k}'-(2k-1)f''s_{f,k}\Big)+k\,\frac{s_{f,k}f^{k}}{(f')^{2k-1}}
	=t_{k+1}+k\,t_k,
	\]
	where the last equality uses Eqn. (\ref{rr}). Both sides of Eqn.~\eqref{eq:trec} are meromorphic, so Eqn. \eqref{eq:trec} holds everywhere. Our aim is to find the order of a pole of $t_k$ at a zero of $f'$ where $f\ne0$. Considering $a\in\mathbb{C}$ with $f(a)\neq0$ and $f'(a)=0$, let $m\ge1$ be the order of the zero of $f'$ at $a$.
	Write $f'=(z-a)^{m}w$ with $w$ holomorphic near $a$ and $w(a)\ne0$. Then
	\[
	\psi=(z-a)^{-m}\,\omega_0(z),\qquad \omega_0:=f/w,\qquad \omega_0(a)=\frac{f(a)}{w(a)}\neq0 .
	\]
	We claim that for every $k\ge1$ there exist a function $\omega_k$, holomorphic near $a$ with $\omega_k(a)=\beta_k\neq0$, such that
	\begin{equation}\label{eq:pole}
	t_k=(z-a)^{-s_k}\,\omega_k,\quad \mbox{~where~} s_k=km+k-1 .
	\end{equation}
	For $k=1$ this is the formula for $\psi$ ($s_1=m$, $\omega_1=\omega_0$).
	Suppose Eqn. \eqref{eq:pole} holds for some $k\ge1$; note $s_k\ge1$. Then
	\[
	t_k'=(z-a)^{-s_k-1}\big(-s_k\omega_k+(z-a)\omega_k'\big),
	\]
	and by Eqn.~\eqref{eq:trec}
	\[
	t_{k+1}=(z-a)^{-(m+s_k+1)}\Big[\omega_0\big(-s_k\omega_k+(z-a)\omega_k'\big)-k\,(z-a)^{m+1}\omega_k\Big].
	\]
	The function in bracket is holomorphic near $a$ and its value at $a$ is $-\omega_0\,s_k\beta_k\neq0$, since $m+1\ge1$. Hence Eqn.~\eqref{eq:pole} holds for $k+1$ with
	$s_{k+1}=m+s_k+1=(k+1)m+(k+1)-1$ and $\beta_{k+1}=-\omega_0\,s_k\beta_k$. This proves the claim.
	Thus, by Eqn.~\eqref{eq:pole}, $s_1<s_2<\dots<s_{n-1}$ and each $t_k$ has a pole of order exactly $s_k$ at $a$.
	In Eqn.~\eqref{eq:S} the term $k=n-1$ has the strictly largest pole order $s_{n-1}$ (note $n-1\ge1$). So the other terms cannot cancel its leading part, and $S_{f,n}-z$ has a pole of order $s_{n-1}=(n-1)m+n-2\ge1$ at $a$. Hence $a$ is a pole of $S_{f,n}$.
	Since $S_{f,n}$ is rational, it has finitely many poles. By above, only finitely many zeros $a$ of $f'$ satisfy $f(a)\neq0$.
	The zeros of $f'$ with $f(a)=0$ are among the zeros of $f=pe^{q}$, i.e.\ among the zeros of the polynomial $p\not\equiv0$, which are finitely many.
	Therefore $f'$ has finitely many zeros.
	Clearly, $f'=(p'+pq')e^{q}=Fe^{q}$ where $F=p'+pq'$ is an entire function. Since $f'\not\equiv0$, the function $F\not\equiv0$, and $F$ has finitely many zeros.
	Let $p_1$ be a polynomial whose zeros, counted with multiplicity, are exactly those of $F$ (take $p_1=1$ if $F$ has no zeros). Then $F/p_1$ is entire and zero-free, so $F/p_1=e^{t}$ for some entire $t$, because the plane is simply connected.
	Therefore $f'=p_1e^{q+t}$, and the lemma holds considering $q_1=q+t$.
\end{proof}

 The next lemma provides certain growth estimates that will be needed later.
\begin{lem}
Let $f(z)=p(z)e^{q(z)}$, where $p$ is polynomial and $q$ is a non-constant entire function. Then $T(r,P_{1,k})=o(T(r,e^q))$ and $T(r,P_{2,k})=o(T(r,e^q))$ outside a set of finite linear measure.	
\label{cha1}
\end{lem}
\begin{proof}
We prove the statement using the principle of mathematical induction. It follows from Eqn.~(\ref{rr}) that the statement is trivially true for $k=1$. In this case, $P_{1,1}(z)=1$ and $P_{2,1}(z)=0$. Thus, $$T(r,P_{1,1})=T(r,P_{2,1})=\mathcal{O}(1)=o(T(r,e^q)).$$ For $k=2$, we have $$P_{1,2}=-(p''+2p'q'+pq''+pq'^2)~~\text{and} ~~ P_{2,2}=P_{1,2}'+P_{1,2}q'.$$ Now,
\begin{align*}
&T(r,P_{1,2})=T(r,-(p''+2p'q'+pq''+pq'^2) )\\
&\leq T(r,(p''+2p'q'+pq''+pq'^2))+\mathcal{O}(1)
\\&\leq T(r,p'')+T(r,p')+T(r,q')+T(r,p)+T(r,q'')+T(r,p)+T(r,q'^2)+ \ln 4+\mathcal{O}(1)
\\&\leq 2T(r,p)+T(r,p')+T(r,p'')+T(r,q')+T(r,q'')+T(r,q'^2)+\mathcal{O}(1)
\\&\leq M \log r + T(r,q')+T(r,q'')+T(r,q'^2)+\mathcal{O}(1)\\
&\qquad\qquad (\mbox{for some $M>0$, as $p$ is a polynomial})
\\& \leq M \log r+o(T(r,e^q))+o(T(r,e^q))+2o(T(r,e^q))+\mathcal{O}(\log r)+\mathcal{O}(1)\\
&\qquad\qquad (\mbox{by Lemma \ref{lem2}})
\\& \leq o(T(r,e^q))\\
& \text{and}\\
&T(r,P_{2,2})=T(r,P_{1,2}'+P_{1,2}q')\\
& \leq T(r,P_{1,2}')+T(r,P_{1,2}q')+\ln 2
\\& \leq T(r,P_{1,2})+\mathcal{O}(\log r+\log T(r,P_{1,2}))+T(r,P_{1,2})+T(r,q')+\ln 2
\\& \leq o(T(r,e^q))+o(T(r,e^q))+o(T(r,e^q))+o(T(r,e^q))\\
&\qquad\qquad (\mbox{by Lemma \ref{lem2} and the fact that  $\lim\limits_{r\to\infty}\frac{T(r,P_{1,2})}{\log r}=\infty$})
\\& \leq o(T(r,e^q)).
\end{align*}
Hence, the statement is true for $k=2$. Assuming that the statement is true for $k=m$, we have $T(r,P_{1,m})=o(T(r,e^q))$ and $T(r,P_{2,m})=o(T(r,e^q))$ outside a set of finite linear measure. It follows from the proof of Lemma \ref{sf1} that $$P_{1,m+1}=(p'+pq')P_{2,m}-(2m-1)(p''+2p'q'+pq''+pq'^2)P_{1,m}$$ and $$P_{2,m+1}=P'_{1,m+1}+mP_{1,m+1}q'.$$ As $-P_{1,2}=p''+2p'q'+pq''+pq'^2$, we have $$P_{1,m+1}=(p'+pq')P_{2,m}+(2m-1)P_{1,2}P_{1,m}.$$ Now
\begin{align*}
T(r,P_{1,m+1})&=T(r,(p'+pq')P_{2,m}+(2m-1)P_{1,2}P_{1,m} )\\
	&\leq T(r,p'+pq')+T(r,P_{2,m})+T(r, (2m-1)P_{1,2}P_{1,m} )+\ln 2
	\\&\leq T(r,p')+T(r,p)+T(r,q')+T(r,P_{2,m})+\\
	&~~~~~T(r,(2m-1)P_{1,2})+T(r,P_{1,m})+2\ln 2
	\\&\leq T(r,p')+T(r,p)+T(r,q')+T(r,P_{2,m})+\\
	&~~~~~T(r,P_{1,2})+T(r,P_{1,m})+2\ln 2+\ln (2m-1)+\mathcal{O}(1)
	\\& \leq M \log r+4o(T(r,e^q))+\mathcal{O}(1)\\
	&\qquad\qquad (\mbox{for some natural number $M$})
	\\& \leq o(T(r,e^q))		
\\ \text{and}\\
T(r,P_{2,m+1})&=T(r,P_{1,m+1}'+mP_{1,m+1}q')\\
	& \leq T(r,P_{1,m+1}')+T(r,mP_{1,m+1}q')+\ln 2
	\\& \leq T(r,P_{1,m+1})+\mathcal{O}(\log r+\log T(r,P_{1,m+1}))+\\
	 &~~~~T(r,P_{1,m+1})+T(r,q')+\ln 2+\ln m+\mathcal{O}(1)
	\\& \leq 5o(T(r,e^q)) \\
	&\qquad\qquad \mbox{(by Lemma \ref{lem2})}
	\\& = o(T(r,e^q)),
	\end{align*}	
showing that the statement is true for $k=m+1$ whenever it is true for $k=m$.
\end{proof}
Note that, if $T(r,f)=o( T(r,e^g))$, where $f$ is an entire function and $g$ is a non-constant entire function, then $$T(r,f')=o (T(r,e^g)).$$
This follows directly from the inequality $$T(r,f')\leq T(r,f)+\mathcal{O}(\log r+\log(T(r,f))).$$ 
The following lemma considers a special case of Lemma \ref{cha1} when $f$ and $f'$ are of the same form. This lemma will be used frequently throughout the paper.
\begin{lem}
Let $f(z)=p(z)e^{q(z)}$ and $f^{\prime}(z)=g(z) e^{h(z)}$, where $p\,(\not\equiv 0)$ and $g\,(\not\equiv 0)$ are polynomials, the functions $q$ and $h$ are non-constant entire such that $t(z)=h(z)-q(z)$ is non-constant. Then $p^kP_{1,k}$ and $p^kP_{2,k}$ can be written as $$p^kP_{1,k}=p(A_{0,k}+A_{1,k}e^t+A_{2,k}e^{2t}+\cdots +A_{2k-2,k}e^{(2k-2)t})$$ and $$p^kP_{2,k}=B_{0,k}+B_{1,k}e^t+B_{2,k}e^{2t}+\cdots +B_{2k-1,k}e^{(2k-1)t},$$ respectively, such that  $T(r,pA_{i,k})=o(T(r,e^t))$  and $T(r,B_{j,k})=o(T(r,e^t))$ outside a set of finite linear measure for $i=0,1,\dots,2k-2$ and $j=0,1,\dots,2k-1$.
	\label{cha2}
\end{lem}
\begin{proof}
According to the above assumption, we have $f(z)=p(z)e^{q(z)}$ and $$f'(z)=(p'(z)+p(z)q'(z))e^{q(z)}=(g(z)e^{t(z)})e^{q(z)},$$ which implies that
\begin{equation}
pq'=ge^t-p'=\phi_1e^t+\phi_1^*,
\label{deri}
\end{equation}
where $\phi_1=g$ and $\phi_1^*=-p'$, and
\begin{equation}
p^2q''=\phi_2e^t+\phi_2^*,
\label{deri2}
\end{equation}  where $\phi_2=pg'+pgt'-p'g$ and $\phi_2^*=p'^2-pp''$.

We prove the statement using the principle of mathematical induction. It follows from Eqn. (\ref{rr}) that the statement is trivially true for $k=1$. In this case, $pP_{1,1}=p$ and $pP_{2,1}=0$ giving that $pA_{0,1}=p$. Similarly, the function $pP_{2,1}(z)=0$ gives that $B_{0,1}(z)=0$. Thus, $T(r,pA_{0,1})=o(T(r,e^t))$ and $T(r,B_{0,1})=o(T(r,e^t))$. For $k=2$, using Eqns. (\ref{deri}) and (\ref{deri2}), we have $$p^2P_{1,2}=-p^2(p''+2p'q'+pq''+pq'^2)=(pp'g-p^2g'-p^2gt')e^t-pg^2e^{2t}.$$ Thus, $A_{0,2}=0$, $A_{1,2}=p'g-pg'-pgt'$, and $A_{2,2}=-g^2$ giving that $T(r,pA_{i,2})=o(T(r,e^t))$ for $i=0,1,~\mbox{and~}2$ (using Lemma \ref{lem2}). Moreover,
\begin{align*}
p^2P_{2,2}&=p^2(P'_{1,2}+P_{1,2}q')\\
&=(p^2P_{1,2})'-2pp'P_{1,2}+(pP_{1,2})(pq')~~(\mbox{as  $(p^2P_{1,2})'=2pp'P_{1,2}+p^2P'_{1,2}$})\\
&=(pA_{1,2}e^t+pA_{2,2}e^{2t})'-2p'(A_{1,2}e^t+A_{2,2}e^{2t})+(A_{1,2}e^t+A_{2,2}e^{2t})(ge^t-p')\\
&=B_{0,2}+B_{1,2}e^t+B_{2,2}e^{2t}+B_{3,2}e^{3t}.
\end{align*}
Clearly, $T(r,B_{j,2})=o(T(r,e^t))$ for $j=0,1,2, \mbox{~and~}3$. Hence, the statement is true for $k=2$. Assuming that the statement is true for $k=m$, we have $$p^mP_{1,m}=pA_{0,m}+pA_{1,m}e^t+pA_{2,m}e^{2t}+\cdots +pA_{2m-2,m}e^{(2m-2)t}$$ and $$p^mP_{2,m}=B_{0,m}+B_{1,m}e^t+B_{2,m}e^{2t}+\cdots +B_{2m-1,m}e^{(2m-1)t},$$ satisfying $T(r,pA_{i,m})=o(T(r,e^t))$ and $T(r,B_{j,m})=o(T(r,e^t))$ for $i=0,1,\dots,2m-2$ and $j=0,1,\dots,2m-1$.
	
It follows from the proof of Lemma \ref{sf1} that
\begin{align*}
	p^{m+1}P_{1,m+1}&=p^{m+1}(p'+pq')P_{2,m}+(2m-1)p^{m+1}P_{1,2}P_{1,m}\\
	&=p(p'+pq')p^mP_{2,m}+(2m-1)pP_{1,2}p^mP_{1,m}\\
	&=p(p'+ge^t-p')(B_{0,m}+B_{1,m}e^t+B_{2,m}e^{2t}+\cdots +B_{2m-1,m}e^{(2m-1)t})\\
	&~~~+(2m-1)p(A_{0,2}+A_{1,2}e^t+A_{2,2}e^{2t})(A_{0,m}+A_{1,m}e^t+\cdots +A_{2m-2,m}e^{(2m-2)t})\\
	&=ge^t(B_{0,m}+B_{1,m}e^t+B_{2,m}e^{2t}+\cdots +B_{2m-1,m}e^{(2m-1)t})\\
	&~~~+(2m-1)p(A_{0,2}+A_{1,2}e^t+A_{2,2}e^{2t})(A_{0,m}+A_{1,m}e^t+\cdots +A_{2m-2,m}e^{(2m-2)t})\\
	&=p((2m-1)A_{0,2}A_{0,m}+(pgB_{0,m}+(2m-1)(A_{0,m}A_{1,2}+A_{1,m}A_{0,2}))e^t+\cdots\\
	&~~~+(pgB_{2m-1,m}+(2m-1)A_{2,2}A_{2m-2,m})e^{2mt})\\
	&=p(A_{0,m+1}+A_{1,m+1}e^t+A_{2,m+1}e^{2t}+\cdots +A_{2m,m+1}e^{2mt})
\end{align*}
Since $T(r,pA_{i,m})=o(T(r,e^t))$ and $T(r,B_{j,m})=o(T(r,e^t))$ for $i=0,1,\dots,2m-2$ and $j=0,1,2,\dots,2m-1$, $T(r,pA_{i,m+1})=o(T(r,e^t))$ outside a set of finite linear measure for $i=0,1,\dots,2m$. Similarly,
\begin{align*}
p^{m+1}P_{2,m+1}&=p^{m+1}(P_{1,m+1}'+mP_{1,m+1}q')\\
&=	(p^{m+1}P_{1,m+1})'-(m+1)p'p^mP_{1,m+1}+mp^mP_{1,m+1}pq'\\
&=(p(A_{0,m+1}+A_{1,m+1}e^t+A_{2,m+1}e^{2t}+\cdots +A_{2m,m+1}e^{2mt}))'\\
&~~~~-(m+1)p'(A_{0,m+1}+A_{1,m+1}e^t+A_{2,m+1}e^{2t}+\cdots +A_{2m,m+1}e^{2mt})\\
&~~~~+m(A_{0,m+1}+A_{1,m+1}e^t+A_{2,m+1}e^{2t}+\cdots +A_{2m,m+1}e^{2mt})(ge^t-p')\\
&=B_{0,m+1}+B_{1,m+1}e^t+B_{2,m+1}e^{2t}+\cdots +B_{2m+1,m+1}e^{(2m+1)t}.
\end{align*}
As $T(r,A_{i,m+1})=o(T(r,e^t))$ and $T(r,A'_{i,m+1})=o(T(r,e^t))$ for $i=0,1,\dots, 2m$, we have $T(r,B_{j,m+1})=o(T(r,e^t))$ outside a set of finite linear measure for $j=0,1,\dots,2m+1$.
Thus, the statement is true for $k=m+1$ whenever it is true for $k=m$. 
\end{proof}
We are now in a position to prove Theorem~\ref{first}.
\begin{proof}[Proof of Theorem~\ref{first}]
\textbf{Sufficiency.} Assume that $f$ is a non-constant entire function of the form $f(z)=p(z)e^{q(z)}$, where $p$ and $q$ are polynomials. Then $$f'(z)=(p'(z)+p(z)q'(z))e^{q(z)},$$ from which we get that $$(f'(z))^{2k-1}=e^{(2k-1)q(z)}(p'(z)+p(z)q'(z))^{2k-1}.$$ Using Eqn. (\ref{sch-def2}), we have
\begin{align*}
S_{f,n}(z)&=z+\sum_{k=1}^{n-1}\frac{(-1)^k}{k!}\frac{s_{f,k}(z)}{(f'(z))^{2k-1}}(f(z))^k\\
	&=z+\sum_{k=1}^{n-1}\frac{(-1)^k}{k!}\frac{e^{(k-1)q(z)}P_{1,k}(z)}{e^{(2k-1)q(z)}(p'(z)+p(z)q'(z))^{2k-1}}(p(z))^ke^{kq(z)} \\
	&\qquad\qquad \mbox{(by Lemma \ref{sf1})}\\
	&=z+\sum_{k=1}^{n-1}\frac{(-1)^k}{k!}\frac{P_{1,k}(z)}{(p'(z)+p(z)q'(z))^{2k-1}}(p(z))^k.
\end{align*}
As $P_{1,k},p' $ and $q'$ are all polynomials, the function $S_{f,n}$ is rational. 

We now check that it is not the identity. If $p$ is non-constant, let $a$ be a zero of $f$ of multiplicity $d$. Put $x=z-a$ and write $f=x^d u$ with $u(a)\ne0$. The recurrence
$$t_1=\frac{f}{f'},\qquad t_{k+1}=\frac{f}{f'}t_k'-k t_k,$$
where $t_k=s_{f,k}f^k/(f')^{2k-1}$, gives
$$t_k(z)=k!\binom{1/d}{k}x+\mathcal{O}(x^2).$$
Consequently,
$$S_{f,n}'(a)=\sum_{k=0}^{n-1}(-1)^k\binom{1/d}{k}
=\prod_{j=1}^{n-1}\left(1-\frac1{dj}\right)\in[0,1),$$
and $S_{f,n}\not\equiv z$. If $p$ is constant, then $q$ is non-constant. If $\deg q\ge2$, choose a zero of $q'$. The local pole-order calculation in the proof of Lemma~\ref{f'finitepole} applies independently of that lemma's rationality hypothesis and shows that $S_{f,n}$ has a pole there. If $q(z)=az+b$, $a\ne0$, a direct calculation gives
$$S_{f,n}(z)=z-\frac1a\sum_{k=1}^{n-1}\frac1k\ne z.$$
Thus $S_{f,n}$ is a non-identity rational function.
	\par
\textbf{Necessity.}
	Conversely, let $S_{f,n}$ be a rational map.  Since all the roots of $f$ are attracting fixed points of $S_{f,n}$, the function $f$ has  at most  finitely many roots. This gives that
	\begin{equation}
	f(z)=p(z) e^{q(z)},\label{1}
	\end{equation}
	where $p$ is a polynomial and $q$ is an entire  function.  We assert that $q$ is also a polynomial.
	\par
	Since $S_{f,n}$ is rational, it has at most finitely many poles. It follows from Lemma \ref{f'finitepole} that
	\begin{equation}
	f^{\prime}(z)=g(z) e^{h(z)}\label{2}
	\end{equation}	
	for some polynomial $g$ and entire function $h$. From Eqns. (\ref{1}) and (\ref{2}), we have
	\begin{equation}\label{q_poly}
	p^{\prime}(z)+p(z) q^{\prime}(z)=g(z)e^{h(z)-q(z)}.
	\end{equation}
	If $q$ is constant, then $f$ is already a polynomial. If $h$ is constant, then $f'=ge^h$ is a polynomial and hence $f$ is a polynomial. If $t=h-q$ is constant, then $p'+pq'=kg$ for some nonzero constant $k$, so $q'$ (being entire and equal to a rational function with polynomial numerator and denominator) is a polynomial. Thus $q$ is a polynomial. We may therefore assume that $q$, $h$, and $t=h-q$ are all non-constant.
	According to the above assumption, we have $f(z)=p(z)e^{q(z)}$ and $$f'(z)=(p'(z)+p(z)q'(z))e^{q(z)}=(g(z)e^{t(z)})e^{q(z)}.$$
	Now
\begin{align*}
S_{f,n}(z)-z&=\sum_{k=1}^{n-1}\frac{(-1)^k}{k!}\frac{s_{f,k}(z)}{(f'(z))^{2k-1}}(f(z))^k\\
	&=\sum_{k=1}^{n-1}\frac{(-1)^k}{k!}\frac{s_{f,k}(z)}{(f'(z))^{k-1}}\left(\frac{f(z)}{f'(z)}\right)^k \\
	&=\sum_{k=1}^{n-1}\frac{(-1)^k}{k!}\frac{P_{1,k}(z)e^{(k-1)q(z)}}{(g(z))^{k-1}e^{(k-1)h(z)}}\left(\frac{p(z)}{g(z)}\right)^ke^{k(q(z)-h(z))}\\ &\qquad \qquad \mbox{(by Lemma \ref{sf1}, and Eqns. (\ref{1}), (\ref{2}))}\\
	&= \frac{1}{(g(z))^{2n-3}} \sum_{k=1}^{n-1}\frac{(-1)^k}{k!}(g(z))^{2n-3-2k+1}P_{1,k}(z)e^{(2k-1)(q(z)-h(z))}(p(z))^k\\
	&= \frac{1}{(g(z))^{2n-3}} \sum_{k=1}^{n-1}\frac{(-1)^k}{k!}(g(z))^{2(n-k-1)}(p(z))^kP_{1,k}(z)e^{(2k-1)(q(z)-h(z))}\\
	&=\frac{1}{(g(z))^{2n-3}} \sum_{k=1}^{n-1}\alpha_kD_k(z),
	\end{align*}
	where $$D_k(z)=(g(z))^{2(n-k-1)}(p(z))^kP_{1,k}(z)e^{(2k-1)(q(z)-h(z))}$$ and $\alpha_k=\frac{(-1)^k}{k!}$.
Assuming  $$D(z)=(g(z))^{2n-3}(S_{f,n}(z)-z),$$ we have	
		\begin{equation}
		 D(z)=\sum_{k=1}^{n-1}\frac{(-1)^k}{k!}D_k(z).
		 \label{D}
		 \end{equation}
Considering $\alpha_i=\frac{(-1)^k}{k!}$ for $i=1,2,\dots,n-1$,  Eqn. (\ref{D}) can be written as
\begin{equation}\label{dif}
\alpha_1D_1(z)+\alpha_2D_2(z)+\cdots +\alpha_{n-1}D_{n-1}(z)-D(z)=0.
\end{equation}
 Thus, we have	
	 \begin{multline}\label{dif1}
	\alpha_1g^{2(n-2)}pP_{1,1}e^{-t}+\alpha_2g^{2(n-3)}p^2P_{1,2}e^{-3t}+\cdots
	 +\alpha_{n-1}g^{2(n-n+1-1)}p^{n-1}P_{1,n-1}e^{-(2n-3)t}-D=0.
	\end{multline}	
Using Lemma \ref{cha2}, we get 	
\begin{multline}\label{dif3}
\alpha_1g^{2(n-2)}pA_{0,1}e^{-t}+\alpha_2g^{2(n-3)}p(A_{0,2}+A_{1,2}e^t+A_{2,2}e^{2t})e^{-3t}+\cdots\\+\alpha_{n-1}g^{2(n-n+1-1)}p(A_{0,n-1}+A_{1,n-1}e^t+A_{2,n-1}e^{2t}+\cdots +A_{2n-4,n-1}e^{(2n-4)t})e^{-(2n-3)t}
-D=0.
\end{multline}	
	We rewrite  the above equality as
	\begin{equation}\label{dif2}
	f_1e^{g_1}+f_2e^{g_2}+\cdots +f_{2n-3}e^{g_{2n-3}}-D=0,
	\end{equation}
	where $g_i=-it$ for each $i=1,2,\dots, 2n-3$. 

The function $f_i=p\sum_{j\geq (\frac{i+1}{2})}^{n-1}\alpha_jg^{2(n-j-1)}A_{2j-i-1,j}$ is the coefficient of $e^{g_i}$ in Eqn. (\ref{dif3}). Note that at least one of $f_i$ is non-zero. Otherwise, $D=0$ gives that $S_{f,n}$ is the identity map.
It follows from Lemma \ref{cha2} that $$T(r,A_{i,j})=o(T(r,e^t))$$
for $i=0,1,\dots, 2n-4$ and for $j=1,2,\dots,n-1$. This gives that $$T(r,f_i)=o(T(r,e^t))$$ for $i=1,2,\dots,2n-3$. Clearly,  $T(r,-D)=o(T(r,e^t))$. Thus, we have $$T(r, f_i)=o(T(r, e^{g_i -g_k}))$$ for all $1 \leq i < k \leq n$. Applying Lemma  \ref{lem1}  to the Equality \eqref{dif}, we have $f_i(z)\equiv0$ for $i=1,2,\dots,n$. In particular $D \equiv 0$.
	Hence, we get $g(z)\equiv0$ (as $S_{f,n}$ is non-identity) implying that the function $f$ is a constant. This is a contradiction and the proof is complete.
	\end{proof}
\subsection{Schr\"{o}der iteration functions of second kind}
We establish another recurrence relation depending upon $K_{f,n}$ itself for a non-constant entire function $f$.
\begin{lem}\label{newr}
If $K_{f,n}$ is Schr\"{o}der's method  of second kind of order $n$, then
\begin{equation}\label{re_kon}
K_{f,n+1}(z)=z-n\frac{z-K_{f,n}(z)}{n-K'_{f,n}(z)}.
\end{equation} 
\end{lem}
\begin{proof}
It follows from Eqn. (\ref{sch2-def11}) that $$z-K_{f,n}(z)=-(n-1)\frac{\left(\frac{1}{f(z)}\right)^{(n-2)}}{\left(\frac{1}{f(z)}\right)^{(n-1)}}.$$ Taking the derivative of $K_{f,n}$ we have  $$K'_{f,n}(z)=n-(n-1)\frac{\left(\frac{1}{f(z)}\right)^{(n)}\left(\frac{1}{f(z)}\right)^{(n-2)}}{\left[\left(\frac{1}{f(z)}\right)^{(n-1)}\right]^2}.$$ Now
	\begin{align*}
z-n\frac{z-K_{f,n}(z)}{n-K'_{f,n}(z)}
&=z-n\frac{-(n-1)\frac{\left(\frac{1}{f(z)}\right)^{(n-2)}}{\left(\frac{1}{f(z)}\right)^{(n-1)}}}
{(n-1)\frac{\left(\frac{1}{f(z)}\right)^{(n)}\left(\frac{1}{f(z)}\right)^{(n-2)}}{\left[\left(\frac{1}{f(z)}\right)^{(n-1)}\right]^2}}\\
&=z+n\frac{\left(\frac{1}{f(z)}\right)^{(n-1)}}{\left(\frac{1}{f(z)}\right)^{(n)}}=K_{f,n+1}(z),
\end{align*}
which completes the desired conclusion.	
\end{proof}

To prove of Theorem~\ref{sec}, the following lemma is needed.

\begin{lem}\label{lem:recursion-leading}
	Let $p\,(\not\equiv 0)$ be a polynomial and $g$ be an entire function. Define
	$u(z)=\frac{e^{-g(z)}}{p(z)}$.	Then, for every integer $k\ge0$, there exists a differential polynomial
	$r_k$ in $g'$, with polynomial coefficients, such that
	\begin{equation}\label{lem:recursion}
	u^{(k)}(z)
	=
	\frac{e^{-g(z)}r_k(z)}{p(z)^{k+1}},
	\end{equation}
	where	$r_0=1$	and
	\begin{equation}\label{lem:r-recursion}
	r_{k+1}
	=
	p\,r_k'-(k+1)p'r_k-g'pr_k.
	\end{equation}
	
	Moreover, if $w=g'$, then for every $k\ge1$,
	\begin{equation}\label{lem:leading}
	r_k=(-pw)^k+L_k,
	\end{equation}
	where $L_k$ is a differential polynomial in
	$w,w',w'',\ldots$, with polynomial coefficients, satisfying $\deg L_k\le k-1$.
	
	Here the differential degree is the total degree in the variables
	$w,w',w'',\ldots$, while the coefficients are regarded as having
	degree zero.
\end{lem}
\begin{proof}
	We first prove Eqn.~\eqref{lem:recursion} by induction on $k$. For $k=0$, we have 	$u=\frac{e^{-g}}{p}
	=
	\frac{e^{-g}r_0}{p},
	~\mbox{and}~
	r_0=1$, so the assertion holds.	Assume that the statement is true for the natural number $k$. This gives that 
	$u^{(k)}
	=
	\frac{e^{-g}r_k}{p^{k+1}}$ for some $k\ge0$. Differentiating gives
	\[
	\begin{aligned}
	u^{(k+1)}
	&=
	\frac{d}{dz}
	\left(
	e^{-g}r_kp^{-(k+1)}
	\right)\\
	&=
	-g'e^{-g}r_kp^{-(k+1)}
	+
	e^{-g}r_k'p^{-(k+1)}
	-
	(k+1)e^{-g}r_kp^{-(k+2)}p'\\
	&=
	\frac{e^{-g}}{p^{k+2}}
	\left(
	-g'pr_k+pr_k'-(k+1)p'r_k
	\right).
	\end{aligned}
	\]
	Therefore, the statement is true for $k+1$ giving that 
	$	u^{(k+1)}
	=
	\frac{e^{-g}r_{k+1}}{p^{k+2}}$, where	$r_{k+1}
	=
	p\,r_k'-(k+1)p'r_k-g'pr_k$. This proves Eqn.~\eqref{lem:recursion} and
	Eqn.~\eqref{lem:r-recursion}.
	
	\medskip
	
	We now prove the leading-term assertion Eqn.~\eqref{lem:leading}. Considering 
	$w=g'$, we claim that	$r_k=(-pw)^k+L_k,~
	\deg L_k\le k-1$ for every $k\ge1$.	For $k=1$, the recursion gives
	\[
	r_1
	=
	p\,r_0'-p'r_0-g'pr_0
	=
	-p'-pw.
	\]
	Thus
	\[
	r_1=(-pw)+L_1,
	~
	L_1=-p',
	\]
	and hence, $\deg L_1=0\le0$. This establishes the initial case. Now suppose that
	\[
	r_k=(-pw)^k+L_k,
	~
	\deg L_k\le k-1
	\]
	for some $k\ge1$. Since $g'=w$, the recursion becomes
	\[
	r_{k+1}
	=
	p\,r_k'-(k+1)p'r_k-pwr_k.
	\]
	Substituting $r_k=(-pw)^k+L_k$	gives
	\[
	\begin{aligned}
	r_{k+1}
	&=
	p\left((-pw)^k+L_k\right)'
	-(k+1)p'\left((-pw)^k+L_k\right)\\
	&\qquad
	-pw\left((-pw)^k+L_k\right).
	\end{aligned}
	\]
	The last term contains $-pw(-pw)^k=(-pw)^{k+1}$. Therefore
	\[
	\begin{aligned}
	r_{k+1}
	&=
	(-pw)^{k+1}\\
	&\quad+
	\Big[
	p\left((-pw)^k\right)'
	-(k+1)p'(-pw)^k
	+pL_k'\\
	&\qquad
	-(k+1)p'L_k
	-pwL_k
	\Big].
	\end{aligned}
	\]
	Define $$
	\begin{aligned}
	L_{k+1}
	:={}&
	p\left((-pw)^k\right)'
	-(k+1)p'(-pw)^k
	+pL_k'\\
	&\quad
	-(k+1)p'L_k
	-pwL_k.
	\end{aligned}
	$$
	We show that $\deg L_{k+1}\le k$. Note that  $\deg\left((-pw)^k\right)=k$. Differentiation does not increase the total differential degree. Indeed,
	if $M=\prod_j\left(w^{(j)}\right)^{\alpha_j}$ has differential degree $\deg M=\sum_j\alpha_j$, then every monomial occurring in $M'$ also has total differential
	degree $\sum_j\alpha_j$. Consequently, $\deg\left(
	p\left((-pw)^k\right)'
	\right)\le k$	and $	\deg\left(
	p'(-pw)^k
	\right)\le k$. By the induction hypothesis, $\deg L_k\le k-1$.	Hence $\deg(pL_k')\le k-1$, $\deg(p'L_k)\le k-1$,	and $\deg(pwL_k)
	\le1+(k-1)=k$.	Thus every term in $L_{k+1}$ has differential degree at most $k$.
	Therefore, $\deg L_{k+1}\le k$. It follows that $r_{k+1}
	=
	(-pw)^{k+1}+L_{k+1},
	~
	\deg L_{k+1}\le k$. The induction is complete. Hence	$r_k=(-pw)^k+L_k,~
	\deg L_k\le k-1$ for every $k\ge1$.
\end{proof}
Suppose that $u^{(k)}\equiv0$ for some $k\ge0$. Then $u$ is a
polynomial. Since $u=1/f$ and $f$ is entire, $u$ has no zeros. A
non-constant polynomial necessarily has a zero, and hence $u$ must be
a non-zero constant. This implies that $f=1/u$ is constant, contrary
to the hypothesis. Therefore
\begin{equation}\label{no_derivative_zero}
u^{(k)}\not\equiv 0,
~ (k\ge0).
\end{equation}
In particular, $K_{f,n}$ is a well-defined meromorphic function.

Now we prove Theorem~\ref{sec}.
\begin{proof}[Proof of Theorem~\ref{sec}]
\textbf{(Sufficiency)} Assume that $f$ is of the form $pe^{q}$ for two polynomials $p(\not\equiv0)$ and $q$. Then $u=\frac{e^{-q}}{p}$. Applying Eqn. \eqref{lem:recursion} with $g=q$, we obtain	$u^{(k)}
		=
		\frac{e^{-q}r_k}{p^{k+1}}$.
		
		Therefore
		\[
		\frac{u^{(n-2)}}{u^{(n-1)}}
		=
		\frac{e^{-q}r_{n-2}/p^{n-1}}
		{e^{-q}r_{n-1}/p^n}
		=
		\frac{p\,r_{n-2}}{r_{n-1}}.
		\]
		Since $p$ and $q$ are polynomials, the recursion Eqn.
		\eqref{lem:r-recursion} shows that every $r_k$ is a polynomial. Moreover,
		by Eqn. \eqref{no_derivative_zero}, we have $r_{n-1}\not\equiv0$. Consequently, $K_{f,n}(z)
		=
		z+(n-1)\frac{p(z)r_{n-2}(z)}{r_{n-1}(z)}$	is rational. This proves the sufficiency.

			\bigskip
		\noindent\textbf{(Necessity)} Assume that for an entire function $f$, 	$K_{f,n}$ is rational. Considering $N=n-1\ge1$, we have 
		\begin{equation}\label{rho_definition}
		\rho=
		\frac{u^{(N-1)}}{u^{(N)}}
		=
		\frac{K_{f,n}-z}{N}
		\end{equation} is also ratioanl. For simplicity, consider $v=u^{(N-1)}$. By Eqn. \eqref{no_derivative_zero}, we have $v\not\equiv0$. Since $v'=u^{(N)}$, Eqn. \eqref{rho_definition} gives $\frac{v'}{v}
		=
		\frac1{\rho}=\sigma$.

		We now prove that $f$ has finitely many zeros. Let $a$ be a zero or a pole of $v$, of order $m\in\mathbb{Z}$.
		Locally,	$	v(z)=(z-a)^m h(z),	~ h(a)\neq0$,	and hence
		\[
		\frac{v'}{v}
		=
		\frac{m}{z-a}+\frac{h'}{h}.
		\]
		Thus every zero or pole of $v$ is a pole of the rational function
		$\sigma$. Therefore $v$ has only finitely many zeros and poles. Now let $a$ be a pole of $u$ of order $s\ge1$. Then $u(z)=(z-a)^{-s}h(z),
		~ h(a)\neq0$.
		
		Differentiating $N-1$ times gives
		\[
		u^{(N-1)}(z)
		=
		C(z-a)^{-s-(N-1)}+\cdots,
		\qquad C\neq0.
		\]
		Thus $u^{(N-1)}=v$ has a pole at $a$. Since $v$ has only finitely
		many poles, $u$ has only finitely many poles. The poles of $u=1/f$ are precisely the zeros of $f$. Hence $f$ has
		only finitely many zeros. Let $p\not\equiv0$ be the polynomial whose zeros, counted with
		multiplicity, are precisely the zeros of $f$. If $f$ has no zeros, we
		take $p=1$. Then $\frac{f}{p}$ is entire and zero-free. Hence there exists an entire function $g$
		such that
		\begin{equation}\label{f_pe_g}
		f=pe^g.
		\end{equation}
Considering $w=g'$, from \eqref{lem:recursion}, we have 
		$
		u^{(N-1)}
		=
		\frac{e^{-g}r_{N-1}}{p^N},
		~\mbox{and}~
		u^{(N)}
		=
		\frac{e^{-g}r_N}{p^{N+1}}.
		$
		Therefore
		\begin{equation}\label{rho_r}
		\rho
		=
		\frac{u^{(N-1)}}{u^{(N)}}
		=
		\frac{p\,r_{N-1}}{r_N}.
		\end{equation}
		
		Writing	$\rho=\frac{a}{b}$, where $a,b$ are non-zero polynomials, from Eqn. \eqref{rho_r},
		\begin{equation}\label{arN}
		a\,r_N=b\,p\,r_{N-1}.
		\end{equation}
		
		By Eqn. \eqref{lem:leading},	$r_N=(-pw)^N+L_N$,
		such that 
		$\deg L_N\le N-1$. Substituting this into \eqref{arN}, we obtain
		\[
		a(-p)^Nw^N+aL_N=bpr_{N-1},
		\]
		and hence,
		\begin{equation}\label{differential_identity}
		a(-p)^Nw^N
		=
		bpr_{N-1}-aL_N.
		\end{equation}
		
		Setting $c:=a(-p)^N$, we have $c\not\equiv0$ is a polynomial, and the right-hand side of
		\eqref{differential_identity} is a differential polynomial in $w$ of
		degree at most $N-1$. Dividing by $c$, we obtain
		\begin{equation}\label{w_equation}
		w^N
		=
		\sum_{\alpha}
		C_\alpha(z)
		\prod_{j\ge0}
		\left(w^{(j)}\right)^{\alpha_j},
		\qquad
		|\alpha|:=\sum_j\alpha_j\le N-1,
		\end{equation}
		where the sum is finite and each $C_\alpha$ is rational. Our aim is to prove that  $T(r,w)=\mathcal{O}(\log r)$. If $w$ is a polynomial, there is nothing to prove. Suppose, therefore, that $w$ is not a polynomial.	For $|z|=r$ and $|w(z)|\ge1$, Eqn. \eqref{w_equation} gives
		$
		|w|^N
		\le
		\sum_\alpha
		|C_\alpha|
		|w|^{|\alpha|}
		\prod_{j\ge1}
		\left|\frac{w^{(j)}}{w}\right|^{\alpha_j}.
		$
		Since $|\alpha|\le N-1$ and $|w|\ge1$, $|w|^{|\alpha|}\le|w|^{N-1}$ and	hence,
		\begin{equation}\label{pointwise_bound}
		|w|
		\le
		\sum_\alpha
		|C_\alpha|
		\prod_{j\ge1}
		\left|\frac{w^{(j)}}{w}\right|^{\alpha_j}.
		\end{equation}
		
		On the other hand, if $|w(z)|<1$, then $\log^+|w(z)|=0$.	Thus, in either case,
		\[
		\log^+|w|
		\le
		\log^+
		\left(
		\sum_\alpha
		|C_\alpha|
		\prod_{j\ge1}
		\left|\frac{w^{(j)}}{w}\right|^{\alpha_j}
		\right).
		\]
		Using $\log^+(x_1+\cdots+x_s)
		\le
		\log s+\sum_{i=1}^s\log^+x_i$	and	$\log^+(xy)\le\log^+x+\log^+y$, and then integrating over $|z|=r$, we obtain
		\begin{equation}\label{proximity_estimate}
		m(r,w)
		\le
		\mathcal{O}(1)
		+
		\sum_\alpha m(r,C_\alpha)
		+
		\sum_{j\ge1}
		A_j
		m\left(r,\frac{w^{(j)}}{w}\right),
		\end{equation}
		where $A_j$ denotes the largest exponent with which
		$w^{(j)}/w$ occurs in Eqn. \eqref{w_equation}.
		
		Since each $C_\alpha$ is rational, $m(r,C_\alpha)=\mathcal{O}(\log r)$.	Consequently, Eqn. \eqref{proximity_estimate} yields
		\begin{equation}\label{proximity_estimate_2}
		m(r,w)
		\le
		\mathcal{O}(\log r)
		+
		\sum_{j\ge1}
		A_j
		m\left(r,\frac{w^{(j)}}{w}\right).
		\end{equation}
		
		By the lemma on the logarithmic derivative, there is a set
		$E\subset(0,\infty)$ of finite linear measure such that, for
		$r\notin E$,	$m\left(r,\frac{w^{(j)}}{w}\right)
		=
		\mathcal{O}\left(\log r+\log^+T(r,w)\right)$. Hence, for $r\notin E$,	$m(r,w)
		=
		\mathcal{O}\left(\log r+\log^+T(r,w)\right)$.	Since $w$ is entire,	$T(r,w)=m(r,w)+\mathcal{O}(1)$. Therefore there exist constants $C,C'>0$ such that
		\begin{equation}\label{T_inequality}
		T(r,w)
		\le
		C\log r+C'\log^+T(r,w),
		\qquad r\notin E.
		\end{equation}
For sufficiently large $T$,	$C'\log T\le\frac12T$.	Thus \eqref{T_inequality} implies $	T(r,w)\le C_1\log r~(r\notin E)$ for some constant $C_1>0$ and all sufficiently large $r$.	Since $E$ has finite linear measure, for every sufficiently large $r$ there exists $r'\in[r,2r]\setminus E$.	As $T(r,w)$ is non-decreasing $T(r,w)\le T(r',w)$. 	Hence	$T(r,w)
		\le
		C_1\log r'
		\le
		C_1\log(2r)
		=
		\mathcal{O}(\log r)$.
		
		Therefore
		\begin{equation}\label{T_log_r}
		T(r,w)=\mathcal{O}(\log r)
		\end{equation}
		for all sufficiently large $r$. Thus, by Eqn. \eqref{max_modulus} and Eqn. \eqref{T_log_r}, we have
		\[
		\log M(r,w)
		\le
		3T(2r,w)+\mathcal{O}(1)
		=
		\mathcal{O}(\log r).
		\]
		This gives that  there exist constants $C>0$ and $r_0>0$ such that
		\begin{equation}\label{polynomial_growth}
		M(r,w)\le Cr^C,
		\qquad(r\ge r_0).
		\end{equation}
		If the Taylor expansion of $w$ is 	$w(z)=\sum_{k=0}^{\infty}a_kz^k$, by Cauchy's estimates,
		$$
		|a_k|
		\le
		\frac{M(r,w)}{r^k}
		\le
		Cr^{C-k}.
		$$ If $k>C$, letting $r\to\infty$ gives	$a_k=0$. Hence only finitely many Taylor coefficients are non-zero, and therefore $w$ is a polynomial. Since $w=g'$	is a polynomial, $g$ is a polynomial. Setting $q:=g$, Eqn.~\eqref{f_pe_g} becomes $f(z)=p(z)e^{q(z)}$. Thus $p$ and $q$ are polynomials, and the proof is complete.
	\end{proof}

\subsection*{Acknowledgements} Gang Liu is supported by Scientific Research Fund of Hunan Provincial Education Department (No. 24A0515).
Soumen Pal acknowledges the support of the Indian Institute of Technology Madras through a Postdoctoral Fellowship (No.F.Acad./IPDF/R12/2025).
\section{Declarations}
\subsection*{Funding} Not applicable.

\subsection*{Conflict of Interests}
The authors declare that they have no conflict of interest regarding the publication of this paper.

\subsection*{Data Availability Statement}
The authors declare that this research is purely theoretical and is not associated with any data.

\subsection*{Code availability} Not applicable.

\end{document}